\documentclass[12pt]{amsart}
 \usepackage{graphicx}
\usepackage{xstring}
\usepackage{forloop}
\usepackage{bbm, amsmath, amsfonts, amscd, latexsym, amsthm, amssymb, graphicx, mathrsfs, comment}
\usepackage[colorlinks=true]{hyperref}

\usepackage{abstract}
\usepackage{mathtools}
\usepackage{caption}
\usepackage{MnSymbol}
\usepackage[dvipsnames]{xcolor}

\usepackage{tikz-cd}
\tikzset{
  symbol/.style={
    draw=none,
    every to/.append style={
      edge node={node [sloped, allow upside down, auto=false]{$#1$}}}
  }
}

\makeatletter
\newif\if@check@engine  \@check@enginetrue 
\makeatother
\usepackage[usenames,dvipsnames]{pstricks}

\newcommand{\nocontentsline}[3]{}
\newcommand{\tocless}[2]{\bgroup\let\addcontentsline=\nocontentsline#1{#2}\egroup}

\usepackage{booktabs}

\newtheorem{theor}{\hspace{1cm}{\sc Theorem}}[section]
\newtheorem{utver}[theor]{\hspace{1cm}{\sc Proposition}}
\newtheorem{sledst}[theor]{\hspace{1cm}{\sc Corollary}}
\newtheorem{lemma}[theor]{\hspace{1cm}{\sc Lemma}}

\newtheorem{problem}[theor]{\hspace{1cm}{\sc Problem}}
\newtheorem*{utver*}{\hspace{1cm}{\sc Proposition}}
\theoremstyle{definition}

\newtheorem{defin}[theor]{\hspace{1cm}{\sc Definition}}
\newtheorem*{defin*}{\hspace{1cm}{\sc Definition}}
\newtheorem{exa}[theor]{\hspace{1cm}{\sc Example}}
\newtheorem{observ}[theor]{\hspace{1cm}{\sc Observation}}
\newtheorem{rem}[theor]{\hspace{1cm}{\sc Remark}}

\newtheorem{quest}[theor]{\hspace{1cm}{\sc Question}}

\newcommand{\sgn}{\mathop{\rm sgn}\nolimits}

\newcommand{\Pic}{\mathop{\rm Pic}\nolimits}

\newcommand{\Vol}{\mathop{\rm Vol}\nolimits}
\newcommand{\Val}{\mathop{\rm Val}\nolimits}

\newcommand{\codim}{\mathop{\rm codim}\nolimits}

\newcommand{\conv}{\mathop{\rm conv}\nolimits}

\newcommand{\Trop}{\mathop{\rm Trop}\nolimits}

\newcounter{idx}

\newcommand{\rotraise}[1]{
  \StrLen{#1}[\slen]
  \forloop[-1]{idx}{\slen}{\value{idx}>0}{
    \StrChar{#1}{\value{idx}}[\crtLetter]
    \IfSubStr{tlQWERTZUIOPLKJHGFDSAYXCVBNM}{\crtLetter}
      {\raisebox{\depth}{\rotatebox{180}{\crtLetter}}}
      {\raisebox{1ex}{\rotatebox{180}{\crtLetter}}}}
}

\renewcommand{\emph}[1]{{\it {\color{NavyBlue} #1}}}

\def\R{\mathbb R}

\def\Z{\mathbb Z}
\def\Q{\mathbb Q}
\def\C{\mathbb C}
\def\CC{({\mathbb C}^\star)}

\def\CP{\mathbb C\mathbb P}
\def\RP{\mathbb R\mathbb P}

\usepackage{xr-hyper}
\makeatletter

\newcommand*{\addFileDependency}[1]{%
\typeout{(#1)}%

\@addtofilelist{#1}
\IfFileExists{#1}{}{\typeout{No file #1.}}
}\makeatother

\begin{document}

\begin{center}{\Large \sc Engineered complete intersections: eliminating variables and understanding topology}

\vspace{3ex}

{\sc Alexander Esterov}
\end{center}

\begin{abstract}
We continue the study of engineered complete intersections (ECI) -- an umbrella generality for a number of important objects in combinatoiral and applied algebraic geometry (such as nondegenerate toric complete intersections, critical loci of their projections, hyperplane arrangements, generalized Calabi--Yau complete intersections, incidence varieties in algebraic optimization, reaction networks).

In this paper, we work on extending to ECIs several classical results about toric complete intersections. This includes symbolic elimination theory, patchworking over $\R$, and computing basic geometric invariants over $\C$.

Our results apply e.g. to eliminating variables in systems of ODEs, such as reaction networks, computing Newton polytopes of discriminants, constructing real polynomial maps and reaction networks with prescribed topology.

Along the way, we assign a cohomology ring to an arbitrary tropical fan, and relate reducible ECIs to arrangements of pairwise intersecting planes in $\R^n$.
\end{abstract}

\vspace{-1ex}

\tableofcontents

\newpage

We work in the complex algebraic torus $\CC^n$ with standard coordinates $z_1,\ldots,z_n$, and identify its characters $z^a:=z^{a_1}\cdots z^{a_n}$ with points of the standard lattice $a\in\Z^n$.
\section{Introduction}
The theory of Newton polytopes studies {\it toric nondegenerate complete intersections} -- varieties in the algebraic torus $\CC^n$ given by several equations 
$$f_1=\cdots=f_k=0,$$
each $f_i$ a general linear combination of a given finite set of characters $$z^{a_{i,1}},\,z^{a_{i,2}},\,\ldots.$$
Geometry of such varieties is determined by their Newton polytopes $N_i:=\conv\{a_{i,1},a_{i,2},\ldots\}\subset\Z^n$. This is important at least for two reasons:

-- At times, a special variety that we encounter, turns out to be a complete intersection of this form. Then its geometric invariant of interest $G$ can be easily read off from Newton polytopes: $$G=g(N_\bullet).$$

-- At times, if we need a special variety with a prescribed geometric invariant $G_0$, we can construct Newton polytopes $N_1,\ldots,N_k$ with the sought value of the invariant $G_0=g(N_\bullet)$. Then the sought variety is defined by generic linear combinations of monomials from $N_1,\ldots,N_k$.

\vspace{1ex}

Of course, applicability of both scenarios is restricted. In \cite{e24b}, we extended the scope of the Newton polytope theory to include e.g. the following important objects:

-- critical loci of projections of ``usual'' complete intersections;

-- hyperplane arrangement complements;

-- generalized Calabi--Yau complete intersections (\cite{gcicy}, \cite{bh});

-- incidence varieties of algebraic optimization;

-- reaction network ODEs.

Our generalization, called an {\it engineered complete intersection} (ECI), is given by several equations 
$$f_1=\cdots=f_k=0,$$
each $f_i$ a fixed linear combination of a given finite set of generic monomials $$c_1z^{a_{1}},c_2z^{a_{2}},\ldots.$$
Note the difference: we have one list of monomials for all equations (the characters possibly repeat, but their coefficients form a generic tuple $c$), and take linear combinations of these generic monomials with arbitrary fixed coefficients $v_{i,j}$ (where $i$ is the number of the equation, and $j$ is the number of the monomial).

It turns out that such ECIs cover the objects listed above, but at the same time still admit a relatively seamless generalization of the Newton polytope theory, with Newton polytopes replaced with so call tropical complete intersections (TCIs).

Such combinatorial data as TCI, besides the exponents $a_i$, should inevitably depend on the coefficient matrix $(v_{i,j})$ of our equations. However, we proved in \cite{e24b} that it depends only on ranks of certain minors of this matrix (i.e. on a relevant matroid).

\begin{rem}
Note that the very possibility to extend the theory of Newton polytopes to ECIs looks like a good luck, because the topology of a general ECI (e.g. of a hypeplane arrangement complement \cite{rybn}) is not determined by its combinatorial data. (This is in contrast to general complete intersections in the classical sense, whose topology is determined y the Newton polytopes.)
\end{rem}

\begin{exa}\label{exa0}
To have a first idea of how it works, let $r, g$ and $b$ be binomials whose Newton polygons are segments of the respective colors on the picture.

\vspace{1ex} 

\begin{center}
\includegraphics[width=4cm]{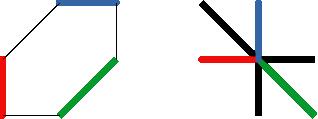}
\end{center}

\vspace{1ex} 

The classical theory of Newton polytopes does not apply to the complete intersection $r+g+b=r+2g+3b=0$ (or to any linear combinations of these equations): for instance, the number of its solutions in $\CC^2$ is 3, twice smaller than predicted by the BKK formula (the lattice area of the Newton polygon on the picture is 6).

However, this complete intersection is engineered, and its geometry can be read off from its TCI. The latter consists of tropical fans $\R^2=:T_0\supset T_1\supset T_2$ and piecewise linear functions $m_1:T_0\to\R$ and $m_2:T_1\to\R$, such that:

-- $m_1$ is the support function of the hexagon;

-- $T_1$ (on the righ) is the corner locus of $m_1$, i.e. the dual fan to the hexagon;

-- $m_2$ equals $m_1$ on the black rays of $T_1$ and equals 0 on the colored rays;

-- $T_2$ is the corner locus of $m_2$ in $T_2$, i.e. the point $\{0\}$ with multiplicity 3.

The multiplicity of $T_2$ reflects the number of points in the complete intersections, while functions $m_i$ (playing the role of Newton polytopes) can be read off from the linear dependencies of the coefficients of the equations.
\end{exa}
In \cite{e24a} and \cite{e24b}, we proved this BKK formula for ECIs, and a number of its extensions.

In this text, we further extend classical methods to ECIs: 

-- develop elimination theory for projections of ECIs (Sections 1 and 2, with applications to discriminants and reaction networks);

-- study characteristic classes and Hirzebruch genera for ECIs (Section 3, which e.g. allows to compute Hodge--Deligne numbers for some SCIs with the Danilov--Khovanskii algorithm);

-- introduce a patchworking-type technique for real ECIs ({\it engineering}, Section 5);

-- relate the question of classifying reducible ECIs to a question of elementary geometry in $\R^n$: how to classify pairs of plane arrangements, such that each plane from the first one intersects each plane from the second one (Section 6).

One tool that we need may be of independent interest: for an arbitrary tropical fan, we introduce a cohomology ring, related to tropical homology of \cite{ikmz} for simplicial fans, and reflecting the actual homology of algebraic varieties for representable fans (so that e.g. the Hodge index theorem and the respective Alexandrov--Fenchel inequality hold; Section 4).

\subsection{Some motivating examples: discriminants and reaction networks}
Given a complete intersection $$f_0=\cdots=f_k=0\eqno{(*)}$$ in the algebraic torus $\CC^n$, its image under the coordinate projection $p:\CC^n\to\CC^k$ is usually a hypersurface. Ignoring its components of higher codimension (if any) and taking the hyperurface components with natural multiplicties, the equation of the resulting Cartier divisor is called the {\it eliminant}, because, from the algebraic perspective, it is the result of eliminating several variables from the system of equations $(*)$. It is much simpler to reason about and compute, if we know its Newton polytope. 

The Newton polytope was computed in \cite{ekh}, if $f_i$'s are generic polynomials with prescribed Newton polytopes (the answer is their {\it mixed fiber body}). This computation was implemented algorithmically in \cite{mm}.
However, the genericity assumptions imposed by \cite{ekh} on $f_i$'s fail in many important examples. 
\begin{exa}\label{exa1}
1. The coordinate projection $\CC^n\to\CC^m$ of the same complete intersection $(*)$ for $m<k$ has the discriminant -- the image of the critical locus $$f_0(z)=\cdots=f_k(z)=0,\,\sum y_i df_i(z)=0\eqno{(**)}$$
(here $y_1,\ldots,y_k$ are new variables, to be eliminated as well). 

\vspace{1ex}

2. Given a system of ODEs 
$$\dot x_i=f_i(x_1,\ldots,x_n),\,i=1,\ldots,n,\eqno{(***)}$$ the elimination of variables in this context aims to find one higher order ODE $x_1^{(N)}=g(x_1^{(N-1)},\ldots,\dot x_1, x_1)$ satisfied by the $x_1$-components of all solutions of the initial functions. One such $g$ (not necessarily minimal in terms of $\deg g$ and $N$) can be found by algebraic elimination of $x^{(j)}_i,\,i>1$, (regarded as independent symbolic variables rather than derivatives) from the initial equations $(***)$ and their symbolic derivatives of all orders up to sufficiently high $N$. Its Newton polytope is estimated from above in \cite{opv} and \cite{mp}.
\end{exa}
Unfortunately, the compete intersections that we have to eliminate in these two examples do not satisfy the genericity assumptions of \cite{ekh}, even if $f_i$'s are generic polynomials with prescribed Newton polytopes.

\vspace{1ex}

The same happens to some other similar questions.

1. How to list the components of complete intersections in Example \ref{exa1}? For discriminants, the question is motivated e.g. by applications to Galois theory (cf.
\cite{compos}). For eliminants of ODEs, it is motivated by the possibility to further simplify the eliminant $g$ of the system $(***)$ (as one of its multipliers is still an eliminant for a general solution of the initial system).

2. If we work over $\R$ instead of $\C$, how to construct complete intersections  with prescribed topology in the setting of Example \ref{exa1}?

Full answers to these questions are given in the literature only for complete intersections that are general in the sense of the classical theory of Newton polytopes, i.e. satisfy roughly the same assumption as in \cite{ekh}: for the classification of components, see \cite{kh16}, and for constructing real complete intersections (in a way classically known for hypersuarfaces as Viro's patchworking), see \cite{sturmfpatchw} and \cite{bihanpatchw}.

\vspace{1ex}

To have a uniform approach to these problems, we have introduced in \cite{e24b} a class of complete intersections (called {\it engineered complete intersections}, ECI), such that

-- the complete intersections of Example \ref{exa1} (as well as many other interesting objects) are themselves coordinate projections of certain nondegenerate ECIs;

-- the classical theory of Newton polytoppes relatively seamlessly extends to nondegenerate ECIs, with so called {\it tropical complete intersections} (TCIs) and {\it matroid complete intersections} (MCIs) playing the role of Newton polytopes.

In this work, we will extend to ECIs classical results of the Newton polytope theory on elimination of variables, enumeration of components, and patchworking. 

In particular, this will allow computing the Newton polytopes of discriminants and eliminants in the general setting of Example \ref{exa1}. 
Enumeration of components and patchworking that we establish for ECIs, do not directly cover Example \ref{exa1} in full generality, since our result require ECI to have some stronger genericity property (being {\it sch\"on}), not satisfied in many instances of Example \ref{exa1}. Nevertheless, this shows that understanding topology in the setting of Example \ref{exa1} essentially boils down to understanding how exactly the complete intersections in Example \ref{exa1} fail to be sch\"on: this question has local nature, and thus a priory simpler than the initial global questions about topology.

\subsection{Engineered complete intersections}

\begin{defin}\label{defeci}
A (symbolic) ECI is a system of equations $$f_i(z):=\sum_{a\in A} v_{a,i} c_a z^a,\eqno{(\diamondsuit)}$$
where $A\subset\Z^n$ is a finite {\it support} multiset, $z^a=z_1^{a_1}\cdots z_n^{a_n}$ is a monomial function on the algebraic torus $\CC^n$, $v_{a,i}$ are complex numbers, and $c_a$ are symbolic coefficients (with respect to given $A$ and $v_{a,i}$).
\end{defin}
\begin{rem}
Instead of indexing monomials with a multiset $A\subset\Z^n$, one can prefer to index them with an abstract finite set mapped (not necessarily injectively) to $\Z^n$, as is done in a related paper \cite{basec}.
\end{rem}
Giving complex values to the symbolic coefficients $c_a$, the systems of equations $f_1=\cdots=f_i=0,\, i=0,\ldots,k$, define nested algebraic sets in $\CC^n$ (with multiplicities), and the tropical fans $T_0\supset\cdots\supset T_k$ of these algebraic cycles are nested. 
If these fans are the same as they would be for almost all other tuples of values of $c_a$, the resulting complete intersection is said to be a {\it nondegenerate ECI}. 

While this nondegeneracy condition may look frightening, it is reduced to explicit conditions on $c_a$'s in \cite{e24a}, Section 4.

\begin{exa}\label{exa2}
1. Consider a generic (or nondegenerate) complete intersection $f_1=\cdot=f_k=0$ in the classical sense, i.e. a tuple of polynomials with prescribed support sets $A_1,\ldots,A_k$ and symbolic (or general in the sense of \cite{kh78} and \cite{ekh}) coefficients.
It is a symbolic (or respectively nondegenerate) ECI with the support multiset $\bigsqcup_i A_i$. 

\vspace{1ex}

2. A hyperplane arrangement complement is an ECI supported at the standard simplex. See \cite{e24a} and \cite{e24b} for many other noteworthy examples.

\vspace{1ex}

3. The right hand side of a reaction network system of ODEs is an ECI.

\vspace{1ex}

4. In the setting of Example \ref{exa1}.1, assume that $$f_1=\cdots=f_k=0\eqno{(*)}$$ is an ECI supported at $A$, and further rewrite each critical locus equation $\sum_i y_i\partial f_i/\partial z_j(z)=0$ as a system $$\sum_i y_i q_{i,j}/z_j=0,\, q_{i,j}-z_j\partial f_i/\partial z_j(z)=0.$$
One may notice that a very similar system $$\sum_i c_{i,j} y_i q_{i,j}/z_j=0,\, q_{i,j}-z_j\partial f_i/\partial z_j(z)=0 \eqno{(\diamondsuit\diamondsuit)}$$
is an ECI: it involves symbolic coefficients of the initial equations $f_i$ and new dummy symbolic coefficients $c_{i,j}$, and has the support set $(A\times\{0\})\cup\{y_i q_{i,j}/z_j, q_{i,j}\}\subset\Z^n\times\Z^N$ (encoding lattice points with the respective monomials). 

One can directly check (using \cite{e24a}, Section 4) that giving sufficiently general values to the coefficients of the initial ECI $(*)$, and value 1 to the dummy coefficients $c_{i,j}$, makes the resulting ECI $(\diamondsuit\diamondsuit)$ nondegenerate. Thus the sought critical locus is the coordinate projection of an ECI. 

\vspace{1ex}

5. The same trick of adding dummy variables equal to logarithmic (higher) derivatives helps in the setting of Example \ref{exa1}.2 as well: for a reaction network ODE system $(**)$ with sufficiently general kinetics (i.e. if $f_1=\cdots=f_n=0$ is an ECI with sufficiently general values of the coefficients), the eliminant is the coordinate projection of an ECI.
\end{exa}

These observations lead us from the motivating problems to the general one.
\begin{problem}\label{probl1} Given a nondegenerate ECI in $\CC^n$:

1. Compute the Newton polytope of its eliminant under the coordinate projection $\CC^n\to\CC^k$;

2. Compute its topological invariants over $\C$, such as characteristic classes and Hirzebruch genera;

3. Understand whether it is irreducible;

4. Learn to patchwork its instances over $\R$, in roder to construct real ECIs with prescribed topology.
\end{problem}

To make this request precise, we should specify in terms of what combinatoiral data of our ECI we want to express the answers.
\begin{defin}
1. A {\it matroid complete intersection} (MCI) is a matroid on a lattice multiset $A\subset\Z^n$.

\vspace{1ex}

2. The ECI $(\diamondsuit)$ is said to represent the MCI $r:2^A\to\Z$, if the matroid $r$ is represented by the vectors $v_{a}=(v_{a,1},\ldots,v_{a,k})$ (i.e., for every $B\subset A$, the vectors $v_{b},\, b\in B$, span a space of dimension $r(B)$).
\end{defin}
Not every MCI is representable (because not every matroid is), but, importantly for us, this relation works well in the other direction: the MCI of a given ECI is easy to read from the settings like our motivating Examples \ref{exa2}.4 and 5.
We will solve Problems 1-4 in terms of MCI of a given ECI in Sections 2, 3, 6 and 5 respectively.

\section{Elimination theory}

\subsection{Tropical complete intersections} The solution of Problem \ref{probl1} regarding elimination consists of two steps:

\vspace{1ex}

i) compute the tropicalization $T$ of an ECI in terms of its MCI;

\vspace{1ex}

ii) compute the eliminant polytope $N$ of a non-degenerate ECI in terms of $T$.

\vspace{1ex}

If the tropicalization is understood in a naive way -- as the tropical fan of the algebraic set defined by the non-degenerate ECI, both steps can be done by references to the existing literature. 
Step (ii) is the theorem of \cite{st} that the sought polytope $N$ is dual to the projection of the tropical fan $T$.

Step (i) is an explicit construction by induction on the number of equations $k$, which we recall in Proposition \ref{proptrop} below.

However, already this construction will require a more clever notion of the tropicalization. It is not only necessary to do step (i) of our solution, but will eventually give a much more practical way to do step (ii) -- a key Corollary \ref{mainth}. 
\begin{defin}[\cite{e24a}]\label{deftrop}
A {\it tropical complete intersection} (TCI) in a tropical fan $T$ is a sequnce of nested tropical fans $T=:T_0\supset T_1\supset\ldots\supset T_k,\, \codim T_i=i$, and piecewise linear functions $m_i:R^n\to\R$, such that each $T_i$ equals the corner locus of the restriction $m_i|_{T_{i-1}}$ (denoted by $\delta(m_i\cdot T_{i-1})$).
\end{defin}
\begin{rem} 1. By saying that a function $m$ is piecewise linear we imply that $m(tx)=tm(x)$ for every $t>0$. 

2. Note that the functions $m_i$ define the sequence of fans $T_i$ uniquely, but not vice versa. Moreover, note that the tail $T_i,\ldots,T_k$ does not depend on values of $m_i$ outside the preceding fan $T_{i-1}$. 

3. In this text, we mostly deal with TCIs in $T=\R^n$, unless otherwise stated.
\end{rem}
\begin{utver}[\cite{e24b}, Section 4]\label{proptrop}
Given an MCI $r:2^A\to\Z$ on a multiset $A\subset\Z^n\subset\R^n$, define its  TCI in $(\R^n)^*$ by induction on the number of fans $k$: assume by induction that the fans $T_i$ and functions $m_i$ are already defined for $i<k$, and, for every primitive covector $l\in T_{k-1}$, define $m_k(l)$ as the maximal $m\leqslant m_{k-1}(l)$ such that the set $\{a\in A\,|\,l(a)\geqslant m\}$ is linearly independent in the sense of the matroid $\min(r,k)$. (This in turn defines the fan $T_k=\delta(m_k\cdot T_{k-1})$.)

If the MCI r is represented by an ECI $f_1=\cdots=f_q=0$, then the fan $T_k$ constructed above is the tropical fan of the set defined by the equations $\tilde f_1=\cdots=\tilde f_k=0$, where $\tilde f_i$ is a general linear combination of the polynomials $f_1,\ldots,f_q$ specialized to general complex values of their symbolic coefficients. In particular, for $k=q$, we get the sought tropical fan of a nondegenerate specialization of the given ECI $f_1=\cdots=f_q=0$.
\end{utver}

\begin{rem}
1. See \cite[Section 7]{e24b} for a more direct application of passing from $f$ to $\tilde f$; note that Section 4 therein gives a slightly more involved inductive procedure of computing the fan $T_q$ with the advantage that each intermediate fan $T_k$ is the tropicalization of $f_1=\cdots=f_k=0$ itself, without tildes. Passing from $f$ to $\tilde f$ make some steps of the procedure tautological, specializing to Proposition \ref{proptrop}.

2. We see that computing the TCI of an ECI from its MCI, besides the sought tropical fan of the ECI, gives for free some extra information: piecewise linear functions $m_i$ on the fan. The meaning of $m_i$ is the tropicalization of the restriction of the function $\tilde f_i$ onto the set $\tilde f_1=\cdots=\tilde f_{i-1}=0$ (this will be made precise below).

3. Not every polynomial restricted to an algebraic set admits a tropicalization in this sense. Those who do, define complete intersections that we call Newtonian. It is a vast generalization of nondegenerate ECIs. It is in this generality that we prove the aforemantioned Corollary \ref{mainth}.
\end{rem}

Before switching to this generality, we discuss how the solution of our elimination problem \ref{probl1} applies to discriminants and resultants of ECIs.

Regard the space $\C^A$ as the space of specializations of our symbolic ECI $(\diamondsuit)$: each point $x\in\C^A$ corresponds to a specialization of $(\diamondsuit)$ with coefficients $c_a:=x_a,\, a\in A$. In this space, nondegenerate ECIs form a Zariski open subset $U$.
\begin{defin}
1. The complement $B:=\C^A\setminus U$ is called the {\it bifurcation set} of the ECI.

2. The closure $D\subset\C^A$ of all specializations having a degenerate solution ($z$ such that $f_1(z)=\cdots=f_k(z)=0$ and $df_1(z)\wedge\cdots\wedge df_k(z)=0$) is called the {\it discriminant} of the ECI.

3. For $k=n+1$, both sets coincide and are called the {\it resultant} of the (-1)-dimensional ECI.
\end{defin}
Clearly, $R$ is the eliminant of the {\it tautological complete intersection} $S$, given in $\C^A\times\CC^n$ by the equations $(\diamondsuit)$ in which both $z$ and $c$ are treated as variables.
This complete intersection is an image of a hyperplane arrangement complement $H$ under an endomorphism of the torus $\CC^{|A|+n}$, thus it is a general ECI (representing the MCI with support equal to $A$ and the matroid defined by $H$). In particular, as we know by now, the Newton polytope of the eliminant $R$ can be computed in terms of this matroid.

From this computation, one can find that the Newton polytope of the resultant $R$ can be expressed in terms of {\it basecondary polytopes} (a version of the secondary polytope twisted with some polymatroid). For the definition and the exact expression for $n=1$, see \cite{basec}.

The same for the discriminant $D$: it is the discriminant of the projection of the tautological general ECI $S$, thus its codimension 1 part can be represented as the eliminant of a nondegenerate ECI in more variables (Example \ref{exa2}.4), thus its Newton polytope can be computed as explained in this section.

\begin{rem}
1. It would be interesting to compute the Newton polytope of the bifurcation set $B$ and to know whether $B\cup D$ always has pure codimension 1. See \cite{adv} for the respective results about nondegenerate complete intersections in the classical sense.

2. It would be interesting to classify MCIs of those ECIs for which the discriminant $D$ is not of pure codimension 1. See \cite{cd} for a survey on this well known problem in the classical setting.
\end{rem}

\subsection{Support functions of eliminant polytopes}

The support function of a polytope $N\subset \R^n$ is a piecewise linear function $N(\cdot):(\R^n)^*\to\R$ sending every linear function $l:\R^n\to\R$ to $\max l|_N$. The polytope may be reconstructed from its support function. So, if we want to compute the eliminant polytope for a Newtonian complete intersection (which includes all eliminants discussed before), it is enough to compute its support function. This will be done in Corollary \ref{mainth}, but first we define NCI.

By the degree of a complex Laurent series $\sum_{k\geqslant d}c_kt^k$ we mean $-\min_{c_k\ne 0} k$, by a germ of an analytic curve $C:\CC\to\CC^n$ a tuple of converging Laurent series $(c_1,\ldots,c_n)$, and set $\deg C:=(\deg c_1,\ldots,\deg c_n)\in\Z^n$.
\begin{defin}[\cite{e24a}]
1. Given an algebraic variety $V\subset\CC^n$, let $|T|\subset\R^n$ be the support of the tropical fan of $V$, i.e. the minimal polyhedral complex containing $\deg C$ for every germ of an analytic curve $C:\CC\to V$.

A regular function $f:V\to\C$ is said to be {\it Newtonian}, if the ratio $\deg f\circ C/\gcd(\deg C)$ depends only on $\deg C$, no matter what curve $C:\CC\to V$ we choose. Assigning this ration to $\deg C$ extends to a unique piecewise linear function on the fan $|T|$: it will be called the tropicalization of $f$.

2. A Newtonian complete intersection (NCI) in the variety $V$ is a sequence $V=V_0\supset V_1\supset\ldots\supset V_k$ and Newtonian functions $f_i:V_{i-1}\to\C$, such that $V_i=\{f_i=0\}$ for every $i$.

3. The tropicalization of this NCI is the TCI given by the tropical fans $T_i$ of $V_i$ and tropicalizations $m_i$ of $f_i$ for all $i$.
\end{defin}
\begin{rem}
1. We will mostly study NCIs in the complex torus $V:=\CC^n$ in this text, unless stated otherwise. Examples of such NCIs include nondegenerate ECIs (and, all the more so, nondegenerate complete intersection in the sense of the classical Newton polytope theory). 

However, there are many more: for instance, a general polynomial $f(x_1,\ldots,x_n)$ gives rise to the symmetric NCI $f(x_1,x_2,x_3,\ldots,x_n)=f(x_2,x_1),x_3,\ldots,x_n)=0$. It naturally appeared in the work of the author and Lionel Lang on Galois theory \cite{el}, and led the author to the subsequent papers on NCIs their various special cases\cite{e24a}, \cite{e24b}, \cite{basec}.

2. Note that, by the definition of TCI (definition \ref{deftrop}), each function $m_i$ completely defines the tropical fan structure of $T_i$ (including the weights) via the expression $T_i=\delta(m_i\cdot T_{i-1})$ (which in turn requires the weights of $T_{i-1}$ by induction on $i$).
\end{rem}
\begin{sledst}[BKK formula for NCIs, \cite{e24a}]\label{bkknci}
In particular, the number of points of a 0-dimensional NCI $f_1=\cdots=f_n=0$ equals the weight of the 0-dimension tropical fan $$ \delta (m_n\cdot \delta( m_{n-1}\cdot \delta(\ldots)))\eqno{(\square)}$$ (supported at $\{0\}$). 
\end{sledst}
\begin{exa}[\cite{e24a}]
If $f_i$'s are general polynomials with prescribed Newton polytopes $N_i$, then, in the respective TCI $(T_i,m_i)$, the functions $m_i$ are the support functions of the polytopes $N_i$. 
\end{exa}
\begin{rem}
1. In particular, in the setting of this example, $(\square)$ expresses the mixed volume $MV(N_1,\ldots,N_n)$ in terms of the support functions of the arguments \cite{bernst}. 

2.Actually, to compute $(\square)$, we do not need individual functions $m_i$, but only the product $m$ of arbitrary piecewise linear extensions of the functions $m_1,\ldots,m_n$ onto the whole $\R^n$: namely $(\square)$ equals
$$\delta^n(m)/n!, \eqno{(\square\square)}$$
where $\delta^n$ is the iterated corner locus of a piecewise polynomial function in the sense of \cite{mmjpoly}.

3. This in particular implies that

-- the expression $(\square)$ is symmetric under any permutation of $m_i$'s (upon extending them arbitrarily onto the whole $\R^n$, so that the permuted expression makes sense);

-- the mixed volume of polytopes with support functions $m_1,\ldots,m_n$ does not depend on values of $m_i$ outside the tropical fan $ \delta (m_{i-1}\cdot \delta( m_{i-2}\cdot \delta(\ldots)))$.
\end{rem}

In order to compute the eliminant polytope of a NCI, we will need the following twisted version of the mixed volume $(\square)$.
\begin{defin}
The {\it mixed shadow volume} of piecewise linear functions $m_0,\ldots,m_n$ in $\R^n_x\times\R_t$ is the number $\delta^{n+1} m/(n+1)!$, where the function $m(x,t)$ is defined as $\prod_i m_i(x,t)-\prod_i m_i(x,0)$ for positive $t$ and $0$ for negative (note that it is continuous if $m_i$ are).
\end{defin}
\begin{defin}
Given a $k$-dimensional tropical fan $T$ in $\R^n\times\R^{k+1}$, define its tropical eliminant as the support function of the polytope in $(\R^{k+1})^*$, dual to the projection of $T$. In other words, the tropical elimiinant is the piecewise linear function in $\R^{k+1}$ whose corner locus equals the projection of $T$.
\end{defin}
\begin{theor}\label{mainthcomb}
Given a TCI $(T_i,m_i),\,i=1,\ldots,n+1$, in $\R^n\times\R^k$, the value of the tropical eliminant at a primitive vector $v\in\R^{k+1}$ is the mixed shadow volume of the restrictions of $m_1,\ldots,m_{n+1}$ onto the space $\R^n\oplus(\R\cdot v)$.
\end{theor}
\begin{sledst}\label{mainth}
If a nondegenerate ECI $f_1=\cdots=f_{n+1}=0$ in $\CC^n\times\CC^{k+1}$ has TCI $(T_i,m_i),\,i=1,\ldots,n+1$, then the support function of its eliminant polytope at a primitive vector $v\in\R^{k+1}$ equals the mixed shadow volume of the restrictions of $m_1,\ldots,m_{n+1}$ onto the space $\R^n\oplus(\R\cdot v)$.
\end{sledst}
\begin{rem}
For a nondegenerate complete intersection in the classical sense, this boils down to the combination of \cite{ekh} (the eliminant polytope is the mixed fiber body) and \cite{mfb} (the formula for the support function of a mixed fiber body). In particular, the latter work defines the mixed shadow volume for convex bodies.
\end{rem}
\begin{proof}
Denote the sought tropical eliminant by $m$, and our candidate for this function (whose values are defined as mixed shadow volumes) by $M$. In order to prove $m=M$, it is enough to prove $\delta(m\cdot R)=\delta(M\cdot R)$ for every 1-dimensional tropical fan $R$. For this, denote $\tilde R:=[\R^n]\times R$, and note that $$\delta(m\cdot R)=\delta(m)\cdot R=\tilde R\cdot T_{n+1}=\tilde R\cdot\delta^n(\prod m_i)/(n+1)!=\delta^{n+1}(\prod m_i(z,t)\cdot \tilde R)/(n+1)!.$$
Here the first and the  last equalities are by the ``Leibnitz rule'' for corner loci (\cite{mmjpoly}), and the second equality is because $\delta(m)$ is the projection is $T_{n+1}$.

To continue this chain of equalities, rewrite it as $\delta^{n+1}([\prod m_i(z,t)-\prod m_i(z,0)]\cdot \tilde R)/(n+1)!$, because $\delta^{n+1}(\prod m_i(z,0))=0$ since the argument depends on less than $n+1$ variables. Finally, rewrite it as $\delta^{n+1}(m\cdot (\tilde R+\tilde R_-))/(n+1)!$, where $\tilde R_-$ is central symmetric to $\tilde R$, and $m$ is $\prod m_i(z,t)-\prod m_i(z,0)$ restricted to $\tilde R$ and then extended by 0 to the remaining part of $\tilde R+\tilde R_-$ (note that the extension is continuous by construction). 

If the rays of $R$ are generated by primitive vectors $v_i$ (repeating to reflect the multiplicities), then $\tilde R+\tilde R_-$ is the sum of the planes $[R^n\times(\R\cdot v_i)]$, so our expression equals the sum of mixed shadow volumes $M(v_i)$, which is the sought number $\delta(M\cdot R)$.
\end{proof}

\section{Refinements of the BKK formula for sch\"on complete intersections} 

In this section, we review some refinements of the BKK formula for NCIs (Corollary \ref{bkknci}), which may help to refine Corollary \ref{mainth} and acquire more combinatorial information about eliminants and discriminants of projected NCIs.

This will require an additional assumption on our NCI.
\begin{defin}
A complete intersection $f_1=\cdots=f_k=0$ in the torus $\CC^n$ is {\it sch\"on} (SCI), if there exists a toric compactification of $\CC^n$, in which the closure of $f_1=\cdots=f_i=0$ is smooth for every $i$ and transversal to every toric orbit it intersects.
\end{defin}

Requiring our complete intersection to be Newtonian and sch\"on at the same time is the same as imposing a purely  algebraical condition (not involving toric compactifications and tropical fans) of being {\it nondegenerte upon cancellations} (NUC).
\begin{defin}\label{defcancel}
1. A cancellation of a tuple of polynomials $(f_1,\ldots,f_k)$ is a tuple $$\tilde f_i(z):=f_i(z)+g_{i,i-1}(z)f_{i-1}(z)+\cdots+g_{i,1}(z)f_1(z).$$

\vspace{1ex}

2. A complete intersection $C\subset\CC^n$ is {\it cancellable} for a linear valuation $v:\Z^n\to\Z$, if there exists a cancellation, whose $v$-leading parts  define a complete intersection. 

If the latter is moreover regular, then $C$ is said to be NUC for $v$, if, moreover.

\vspace{1ex}

3. A complete intersection is cacnellable/NUC, if it is for every non-trivial valuation.
\end{defin}
\begin{theor}[\cite{e24b}, Theorem 6.1]
1. A cancellable complete intersection is NCI.

\vspace{1ex}

2. Being NUC is equivalent to being Newtonian and sch\"on.
\end{theor}

All examples of NCIs that we discussed before (including general ECIs) are at the same time SCIs (and thus NUC).
\begin{rem}
However, nondegeneracy of an ECI does not imply being sch\"on.
\end{rem}

\subsection{CSM classes} We shall compute tropical characteristic classes of a sch\"on NCI $f_1=\cdots=f_k=0$, in the sense of \cite{e13} and \cite{e24c}. These combinatorial entities encode in particular CSM classes of our set (as elements in the homology of any toric compactification), including the Euler charcteristics (represented by the 0-th class).

\begin{theor}[\cite{e24a}]\label{euleci}
If a sch\"on NCI $V:=\{f_1=\cdots=f_k=0\}$ has TCI $(T_i,m_i)$, then its tropical characteristic class equals the expression $$\prod_i \frac{\delta m_i}{(1+\delta m_i)},$$
where $m_i$ are arbitrarily extended to continuous piecewise linear functions on the ambient space $\R^n$, the rational function is symbolically expanded as a power series, and each term in this power series is interpreted as an iterated corner locus, similarly to $(\square)$.
\end{theor}
\begin{proof} The result of the computation in the statement is a formal sum of tropical fans $c_i$ of dimensions form $i=0$ to $\dim V$, and we want to prove they represent the actual tropical characteristic classes $C_i$ of $V$.

(0) For $c_0$ it means that the weight of this 0-dimensional fan equals the Euler characteristics of $V$. 
This equality is proved in \cite{e24a}. The full statement reduces to this case as follows.

(1) The statement is proved fora general complete intersections (in the classical sense, given by general polynomials with prescribed newton polytopes) $U\subset\CC^n$ in \cite{e13}.

(2) Both the actual tropical characteristic classes, and the candidate expressions in the statement of the theorem, multiply when we take the intersection $U\cap V$: for the candidate expression this is by definition, and for actual classes see \cite{e13}.

(3) By the nondegeneracy of the intersecion pairing in the ring of tropical fans, it is enough to prove the theorem for $U\cap B$ with any such $U$, but only for the 0-dimensional class.

(4) This reduces to (0), because $U\cap B$ is itself a sch\"on NCI. 
\end{proof}

\begin{rem}\label{remproj}
Knowing tropical/CSM classes of a complete intersection is useful for the study of the geometry of the eliminant (especially its singularity strata), because these classes behave controllably under projections. Specifically, recall that the MacPherson image of a variety $U\subset N$ under a regular map $p:N \to M$ is a function $M\to\Z$ sending every $v\in M$ to the Euler characteristics of $U\cap p^{-1}(v)$. This function can be represented as a linear combination of indicator functions of varieties $V_i\subset M$. Then the pushforward of the CSM class of $U$ (or, stronger, of its tropical characteristic class in case $p$ is a projection $\CC^n\to\CC^m$) is the respective linear combination of the CSM/tropical classes of $V_i$'s.
\end{rem}

\subsection{Hirzebruch genera} We compute them for a sch\"on NCI in terms of its TCI. The answer requires some vocabulary of lattice polytope valuations.

Convex lattice polytopes in $M\simeq\Z^n$ form a semigroup $\mathcal{P}_+$ with respect to  Minkowski summation $A_1+A_2:=\{a_1+a_2,|\,a_i\in A_i\}$. This semigroup embeds into the group $\mathcal{P}$ of continuous piecewise linear functions $M\to\Z$, once we send a polytope $P$ to its {\it support function} $P(v):=\max v|_P$.
The counting valuation $P\mapsto|P\cap M|$ and volume $P\mapsto n!\Vol P$ extend to degree $n$ polynomials $\#$ and $\Vol_\Z:\mathcal{P}\to\Z$ (in the sense that their restriction to every finite-dimensional subgroup of $\mathcal{P}$ is a polynomial). 
Their values outside $\mathcal{P}_+$ can be computed e.g. from the {\it universal valuation} $P\mapsto 1_P$, where $1_P:M\otimes\R\to\{0,1\}$ is the {\it indicator function}: $1_P^{-1}(1)=P$. 
The universal valuation extends to $\Val:\mathcal{P}\to\mathcal{M}:=\{$finite linear combinations of indicator functions of convex polytopes$\}$, as shown in \cite{khp}. The extension is constructive: 

1. Represent any $m\in\mathcal{P}$ as the difference of support functions of convex polytopes $m_+$ and $m_-$ (this is possible: $\mathcal{P}$ is the Grothendieck group of the semigroup $\mathcal{P}_+$);

2. Define the function $\Val m:M\otimes\R\to\R$ as $\sum_{f}(-1)^{dim\, m_--\,\dim f}1_{m_++f}$ over all faces $f$ of the polytope centrally symmetric to $m_-$. It does not depend on the choice of $m_\pm$; in particular:

$$\Vol_\Z m=\int_{M\otimes\R}\Val m=\sum_{f}(-1)^{dim\, m_--\,\dim f}\Vol_\Z({m_++f}),$$ $$\#m=\sum_{a\in M}\bigl(\Val m\bigr)(a)=\sum_{f}(-1)^{dim\, m_--\,\dim f}\#({m_++f}).$$
\begin{rem}
Alternatively, we could use the lattice version of McMullen's polytope algebra $\mathcal{\hat M}:=\mathcal{M}/($the action by parallel translations of $M)$ as the range for our universal valuation \cite{mcm}. However, we could not regard $\Val m$ as a function on $M\otimes\R$, and would not see the useful interpretation of $\Vol_\Z m$ and $\# m$ in the middle part of the equalities above.\end{rem}

For a vector $\beta$, denote the sum of entries by $\sum\beta$, and the number of non-zero entries by $|\beta|$.
\begin{theor}\label{thgenus}
For a smooth Newtonian SCI $H_k\subset\cdots\subset H_1\subset\CC^n$, extend the tropicalizations of its equations to $m_1,\ldots,m_k\in\mathcal{P}$. Then the genus 
$\chi_p(H_k)$ equals
$$\chi_p(m_1,...,m_k):=(-1)^{n-p} 
\sum\nolimits_{\beta\in\Z^k_{\geqslant 0}}(-1)^{\sum\beta}{n+|\beta|\choose p-\sum\beta-|\beta|}\#\Bigl(\sum_i\beta_i m_i\Bigr)\;\mbox{ for }\;p=0,\ldots,n-k.$$

In particular, the Euler characteristics $e(H_k)$ equals $\sum_{\beta\in\Z^k_{\geqslant 0}}(-1)^{\sum\beta}{n+|\beta|-1\choose k+\sum\beta-1}\Vol_\Z\Bigl(\sum_i\beta_i m_i\Bigr)$.
\end{theor}
\begin{rem} These formulas can be further convoluted to 
$$\#(y-1)^n\prod\nolimits_{i=1}^k\frac{1-[m_i]}{1-y[m_i]}\mbox{ and }\Vol_\Z\prod\nolimits_{i=1}^k\frac{m_i}{1-m_i},$$
where the first expression is evaluated in the ring $\mathcal{\hat M}[[y]]$ as elaborated in \cite{gross}, and the second one as in Theorem \ref{euleci} (with the mixed volume of $n$ elements of $\mathcal{P}$ understood as the polarization of the degree $n$ polynomial $\Vol_\Z:\mathcal{P}\to\Z$).
\end{rem}

\noindent {3.3.\, \bf Proof of Theorem \ref{thgenus}.} 
First, neglect the toric specifics and study any complete intersection $H_k\subset\ldots\subset H_1\subset H_0$ with a compactification $H_0\subset X$ such that

i) $\bar H_{i+1}$ is a smooth hypersurface in $\bar H_i$ (bars denote closure in $X$, indices belong to $[0,k]$);

ii) The set $\bar H_{i}\setminus H_i$ is a normal crossing divisor $D_i$;

iii) The divisor $\bar H_{i}$ (with multiplicity 1) is linearly equivalent in $\bar H_{i-1}$ to a certain divisor supported at $D_{i-1}$ (with some multiplicities of components).

Thus both divisors define the same invertible sheaf on $\bar H_{i-1}$, denote it by $I_{i}$.
(Recall that we assign positive multiplicities to poles, so $I_{i}$ has a section whose divisor is $-\bar H_{i}$.)

Let $\Omega^p_i$ be the sheaf of $p$-forms on $\bar H_i$, whose restriction to every component of $D_i$ is identical 0 (as a $p$-form on that component). We aim at computing $\chi(\Omega^p_k)=\chi_p(H_k)$.
\begin{lemma}[see e.g. \cite{dkh}, (4.1)] For every invertible sheaf $I$ on $\bar H_{i-1}$, we have
$$\chi(\Omega^p_i\otimes I)=\sum\nolimits_{j\geq 0} (-1)^j \Bigl(\chi(\Omega^{p+j+1}_{i-1}\otimes I\otimes I_{i}^{j+1})-\chi(\Omega^{p+j+1}_{i-1}\otimes I\otimes I_{i}^{j})\Bigr).$$
\end{lemma}
Extending each $I_i$ to the whole $X$, we can apply this lemma inductively by $i$.
\begin{sledst}\label{chichi}
$$\chi_p(H_k)=\chi(\Omega^p_k)=\sum_{\beta\in\Z^k_{\geq 0}}(-1)^{\sum\beta}\sum_{\gamma\in\{0,1\}^k}(-1)^{k-\sum \gamma}\chi(\Omega_0^{p+\sum\beta+k}\otimes I_1^{\beta_1+\gamma_1}\otimes\cdots\otimes I_k^{\beta_k+\gamma_k}).$$
\end{sledst}
Back to the setting of Theorem \ref{thgenus}, the SCI $H_k\subset\ldots\subset H_1$ is equipped with:

-- the variety $H_0=T$, an algebraic torus with the character lattice $M\simeq\Z^n$, 

-- its closure $X=X_\Sigma$, a tropical compactification of the SCI $H_\bullet$, 

-- its boundary $D_0$, the union of closures of $(n-1)$-dimensional orbits $D_v\subset X_\Sigma$, indexed by the primitive generators $v$ of 1-dimensional cones in $\Sigma$. 
Finally, a natural surjection
$${\mathcal I}:\mathcal{P}_\Sigma:=(\mbox{piecewise linear functions on the fan }\Sigma)\to\Pic(X_\Sigma),$$
such that a suitable section $s(m)$ of the line bundle $\mathcal{I}(m)$ has the divisor of zeroes and poles $-\sum_v m(v)D_v$, allows to define the line bundle $I_i$ on $X_\Sigma$ as $\mathcal{I}(m_i)$. Recall that the function $m_i\in\mathcal{P}_\Sigma$ is chosen to extend the tropicalizattion of the defining equation $f_i$ of the SCI.
\begin{lemma}
An SCI $H_\bullet$ with $X:=X_\Sigma$ and $I_i:=\mathcal{I}(m_i)$ as above satisfies (i-iii).
\end{lemma}
\begin{proof}
(i) and (ii) are the definition of SCI. The divisor of zeroes and poles of $f_i$ on $\bar H_{i-1}$ equals $-\bar H_i+\sum_{v}m_i(v)D_v\cap\bar H_{i-1}$ by definition of $m_i$. Thus the section $s(m_i)$ of $\mathcal{I}(m_i)$, when restricted to $\bar H_{i-1}$, has the divisor of zeroes and poles linearly equivalent to $-\bar H_i$: this proves (iii).
\end{proof}

To prove Theorem \ref{thgenus}, consider its right hand side as a polynomial $A(m_1,\ldots,m_k)$ on $\mathcal{S}:=\mathcal{P}_\Sigma\oplus\ldots\oplus\mathcal{P}_\Sigma$, and the right hand side of Corollary \ref{chichi} as a polynomial function $B(I_1,\ldots,I_k)$ on $\Pic(X_\sigma)\oplus\cdots\oplus\Pic(X_\sigma)$. Since $B(\mathcal{I}(m_1),\ldots,\mathcal{I}(m_k))=A(m_1,\ldots,m_k)$ on the full dimensional cone $\{$tuples of convex functions$\}\subset\mathcal{S}$ (by \cite{dhn}), these functions equal as polynomials.

\section{"Another tropical homology"}

\subsection{Formalism}
Cohomologically tropical varieties are varieties whose topological cohomology can be understood from their tropical fan. The story started from \cite{zharkov13} proving that hyperplane arrangement complements are CTV, and the term was introduced in \cite{shaw23}.

We present a different point of view on this notion, using the construction of tropical homology given in \cite{mmjpoly}[Section 4]. This will be more convenient for applications to intersection theory in the next sections.

Every tropical fan $T$ in $\R^n$ gives rise to a ring $H(T)$ that has a good reason to be called its tropical cohomology (to be distinguished from tropical homology of smooth tropical varieties introduced by \cite{ikmz}). 
\begin{utver}
Let $B$ be the ring of continuous piecewise polynomial functions on $\R^n$, graded by the degree of polynomials. Let $H^k(T)$ be the image of the map $w_k:B^k\to \{$codimension $k$ tropical fans in $T\}$ defined as $w_k(M):=\delta^k(M\cdot T)/k!$ (where $M\cdot T$ is a tropical fan with piecewise polynomial weights, and $\delta$ stands for the corner locus in the sense of \cite{mmjpoly}). Then the multiplcation in $B$ pushes forward to the multiplication on $$H(T):=\bigoplus_k H^k(T).$$
\end{utver}
\begin{proof}
It is enough to prove that $w(M)=0$ implies $w(MM')=0$ for any $M'\in B^{k'}$. This is because $\delta^{k+k'}(MM'T)/(k+k')!=\delta^{k'}(M')\delta^k(M)T/(k!k'!)=\delta^{k'}(M')\cdot 0$ by the Leibnitz rule for $\delta$ (\cite{mmjpoly}).
\end{proof}
\begin{defin}
The ring $H(T)$ is called the {\it tropical cohomology} of the fan $T$.
\end{defin}
The rationale behind this name is as follows. For a toric fan $\Sigma$, denote the respective toric variety by $X_\Sigma$, the ring of piecewise polynomial functions that are polynomial on each cone of $\Sigma$ by $B_\Sigma$, and its image in $H(T)$ by $H_\Sigma(T)$.
\begin{utver}
Let $T$ be the tropical fan of an algebraic variety $V\subset\CC^n$, and let a smooth toric variety $X_\Sigma\supset\CC^n$ be its tropical compactification (i.e. the closure $\bar V\subset X$ properly intersects all orbits, so that the codimension of the intersection is the sum of the codimensions).

Then $H^\bullet_\Sigma(T)$ is isomorphic to the image of $H^\bullet(X)$ in $H^\bullet(\bar V)$ (induced by the inclusion $\bar V\subset X$).
\end{utver}
\begin{proof}
Consider the diagram 
\begin{tikzcd}
B_\Sigma \arrow[r, "b"] \arrow[d, "w"] & H(X) \arrow[d, "j^*"] \\
H_\Sigma(T)                            & H(\bar V)            
\end{tikzcd}

Here $j^*$ is induced by $j:\bar V\hookrightarrow X$, $w$ is defined above, and Brion's surjection $b$ can be characterized by the property that the intersection number of $b(M)$ with the closure of every $k$-dimensional orbit $O\subset X$ equals the weight of the $O$-cone in the tropical fan $\delta^k(M)/k!$ (where $k:=\deg M$).

We see that $b(M)\in\ker j^*$ iff $b(M)\cap[\bar V]=0$ iff $(b(M)\cap[\bar V])\cdot\bar O=0$ for every orbit $O\subset X$. The latter equality, by our characterization of Brion's surjection, means $\delta^k(M)\cdot T=0$. By the Leibnitz rule for $\delta$, this is equivalent to $\delta^k(M\cdot T)=0$, i.e. $M\in\ker w$.
\end{proof}
\begin{rem}
The map $w$ factors through $B\to B(T)\to H(T)$, where $B(T)$ is the ring of continuous piecewise linear functions on the support of the tropical fan $T$. In some cases, it is more convenient to work with.
\end{rem}
\begin{defin}
A subvariety $V\subset\CC^n$ is said to be {\it homologically toric} (HTV), if the map $H(X)\to H(\bar V)$ is surjective for a given (or, equivalently in the case of sch\"on varieties, any) smooth tropical compactification $X\supset\CC^n$.
\end{defin}
\begin{exa}
1. A curve is HTV iff it is rational.

2. A hyperplane arrangement complement, i.e. an affine plane in $\CC^n$, is HTV.
\end{exa}
\begin{rem}
It would be interesting to classify HTVs in small dimension (especially 2): as we shall see now, they are as convenient as affine spaces in the role of patchworking pieces, if we patchwork a variety to study its intersection theory. 
\end{rem}
\begin{sledst}
1. For an HTV $V$ with the tropical fan $T$ and tropical compactification $\bar V\subset X_\Sigma$, the tropical cohomology $H_\Sigma(T)$ equals the topological cohomology $H(\bar V)$.
The isomorphism sends every cohomology class $\alpha$ to the tropical fan, whose weight on a cone $C\in\Sigma$ equals the value of $\alpha$ on the fundamental class of $\bar V\cap\overline{C-\mbox{\rm orbit}}$.

2. In particular, for any two subvarieties $U_1$ and $U_2$ in $V$, the product of their fundamental classes in $H(\bar V)$ can be recovered from their tropical fans (and the ambient tropical fan of $V$).
\end{sledst}
\begin{proof}
(1) is obvious from the preceding definition and proposition. For (2), represent the tropical fan of $U_i$ as $\delta^{k_i}(M_i\cdot T)/k_i!,\, k_i:=\dim U_i$, then the sought product is defined by the tropical fan $\delta^{k_1+k_2}(M_1\cdot M_2\cdot T)/(k_1+k_2)!$.
\end{proof}
\begin{rem}\label{tropofcycle}
The map $\Trop:H(\bar V)\to\{$tropical fans in $T\}$, described in this corollary, exists even if $V$ is not HTV. However, in that case, the map does not have to be injective, and its image may exceed the tropical cohomology $H_\Sigma(T)$.

An example is given by the fundamental cycles of two intersecting lines in a general quadric $V\subset\CC^3\supset\CP^3$: their tropicalization is the same  (the tropical line), and is not in the tropical cohomology of $T=2\cdot($tropical plane$)$ with integer coefficients.

As we see, omitting the HTV assumption in the preceding corollary, $\Trop U_i$ may miss $H(T)$ at all (so that the intersection number $\Trop U_1\cdot\Trop U_2$ is not defined), or belong to $H(T)$, though the well defined $\Trop U_1\cdot\Trop U_2$ does not reflect the intersection number of the varieties $U_i$. (An example of the latter is given by each $U_i$ being a pair of non-intersecting lines in the quadric.)
\end{rem}
\subsection{Applications} We now give applications to complex and tropical intersection theory, with a view towards ECIs.
\begin{defin}
1. A {\it twist} of a toric variety $X_\Sigma\supset\CC^n$ is a toric variety $\Sigma_{j^{-1}\Sigma}$ for any embedding of lattices $j:\Z^n\to\Z^n$. Note that the twist comes with a map to $X_\Sigma$, induced by $j$.

2. Given a HTV $V\subset\CC^n$ and cycles $\alpha_i$ in the homology of a tropical compactification $\bar V\subset X$, for any twist $j:Y\to X$, the tuple $(Y,j^{-1}\bar V,j^*\alpha_\bullet)$ is called a {\it twisted HTV}.
\end{defin}
\begin{rem}\label{remthtv}
1. Twisted HTV is in general not an HTV (by far): for instance, the surface $\sum_{i=0}^n x_i^d=0$ in $\CP^n$ is a twisted HTV.

2. Still, for a twisted HTV $(Y,W,\alpha_\bullet)$, we have the equality 
$\circ_i\Trop\alpha_i=\circ_i\alpha_i$ (induced from the underlying HTV).
\end{rem}
\begin{defin}
1. Given a variety $V\subset\CC^n_x$, its tropical compactification $X_\Sigma$, and cycles $\alpha_i\in H(\bar V)$, a {\it deformation} of this triple is a variety $V'\subset \CC^n_x\times\CC^1_t$, its tropical compactification $X_{\Sigma'}$, and cycles $\alpha'_i\in H(\bar V')$, such that:

-- $X_\Sigma=X_{\Sigma'}\cap\overline{\{t=1\}}$,

-- $V=V'\cap\{t=1\}$,

-- $\alpha_i=\alpha'_i|_{\bar V}$.

\vspace{1ex}

2. The deformation is said to be {\it cohomologically trivial}, if, for every 1-dimensional cone $r\in \Sigma'$, lying in the open upper halfspace $\R^n\times\R^1_{>0}$, the following forms a twisted HTV:

-- the closure of the $r$-{orbit} $O\subset X_{\Sigma'}$,

-- its intersection with $\bar V'$,

-- the restrictions of $\alpha'_i$ onto this intersection.
\end{defin}
\begin{rem}
A cycle $\alpha$ has a deformation only if its tropicalization (Remark \ref{tropofcycle}) is in the tropical cohomology.
However, we recall from  Remark \ref{tropofcycle} that having $\Trop\alpha$ and $\Trop\beta$ in the tropical cohomology does not assure that one can reconstruct the intersection number $\alpha\cdot\beta$ from the fans $\Trop\alpha$ and $\Trop\beta$.
\end{rem}

\begin{theor}\label{thdef}
Let $V\subset\CC^n$ be any variety with a tropical compactification $\bar V\subset X_\Sigma$. 
If cycles $\alpha_i\in H(\bar V)$ have a common deformation that is cohomologically tropical, then their intersection number equals the intersection number of the fans $\Trop\alpha_i$ in the fan $T:=\Trop V$. (Recall that the latter evaluates as $\delta^{\sum k_i}(\prod_i M_i\cdot T)/(\sum k_i)!$ upon representing $\Trop \alpha_i$ as $\delta^{k_i}(M_i\cdot T)/k_i!$.)
\end{theor}
\begin{proof}
We want to compute $[\bar V']\cdot[\overline{\{t=1\}}]\cdot\prod \alpha_i$. Sending $t\to\infty$, we can replace $[\overline{\{t=1\}}]$ with a linear combination of $[\bar O_i]$, where $O_i$ are the codimension 1 orbits of the deformation toric variety $X_{\Sigma'}$, corresponding to the rays of $\Sigma'$ in the upper halfspace $\R^n\times\R^1_{>0}$. Then the intersection number splits into $[\bar V']\cdot[\bar O_i]\cdot\prod \alpha_i$, which equal the respective tropical intersection numbers by Remark \ref{remthtv}.2.
\end{proof}

\begin{theor}\label{thalf}
Let $T$ be the tropicalization of an irreducible surface $S\subset\CC^n$. 
Then the tropical cohomology $H(T)$ satisfies

1) Hodge index theorem: the intersection pairing on $H(T)$ has signature $(+\,-\,-\,\ldots)$;

2) Alexandrov--Fenchel inequality: for any cycles $\alpha_1$ and $\alpha_2\in H(T)$, if $\alpha_1^2\geqslant 0$, then $(\alpha_1\cdot \alpha_2)^2\geqslant\alpha_1^2\alpha_2^2$.
\end{theor}
\begin{rem}
1. If $F$ is the intersection product of the dual fans of polytopes $\Delta_3,\ldots,\Delta_{n}$, and $\alpha_i$ is the support function of a polytope $\Delta_i,\,i=1,2$, then we get the classical Alexandrov--Fenchel inequality for mixed volumes of $\Delta_1,\ldots,\Delta_{n}$, hence the name. 

2. If a fan is not a tropicalization of an irreducible surface, the theorem may fail. For instance, let $M$ and $N$ be the dual fans of two transversal polygons in $\R^4$, and let $F=M^2+N^2$ be the union of two transversal planes with multiplicities $m$ and $n$ equal to the areas of the polygons. Then $(MNF)^2=M^3N+MN^3=0<m^2n^2=(MMF)^2(NNF)^2$.
\end{rem}
\begin{quest}
Does the theorem, or at least its second part with a stronger assumption on $\alpha_i$, hold true for any 2-dimensional tropical fan $F$ in $\R^n$, such that $F\setminus\{0\}$ is connected?

Note that the question reduces to the case $n=4$ (by a general projection $\R^n\to\R^4$, since its restriction to $F$ is injective).

On the other hand, note that the answer is positive for $n<4$ (because every fan in $\R^3$ is the tropicalization to a surface: the sought surface is the zero locus of a general polynomial with the Newton polytope dual to $F$).
\end{quest}
The proof of the theorem is an update of the classical Teissier--Khovanskii proof of the Alexandrov--Fenchel inequality for mixed volumes.
\begin{proof}
i. As (1) and (2) are equivalent, we shall prove (2). Choose a tropical compactification $\bar S\subset X_\Sigma$, and let $p:\tilde S\to\bar S$ be its resolution of singularities. The resolution has 1-dimensional fibers over finitely many points $R\subset\bar S$.

ii. Represent the tropical fan $\alpha_i$ as the tropicalization of some divisor on $X_\Sigma$, restricted to $\bar S$. This is done as follows: represent it as
$\delta(m_i\cdot T)$ for a continuous piecewise linear $m_i:\R^n\to\R$, represent $m_i$ as the difference of support functions of polytopes $P_i$ and $Q_i$, let $f_i$ be the ratio of general polynomials with the Newton polytopes polytopes $P_i$ and $Q_i$. 
Then $f_i$ defines a divisor in $\CC^n$, and its closure in $X_\Sigma$ is the sought divisor denoted by $H_i$.

By construction, it has the following properties:

(1) $\alpha_i=\delta(m_i\cdot T)=\delta(\Trop f_i\cdot T)$;

(2) $\alpha_i\cdot \alpha_j=\delta^2(\Trop f_i\cdot \Trop f_j \cdot T)=($the intersection number of $\bar S, H_i$ and $H_j$ in $X_\Sigma)=($the intersection number of $H_i|_{\bar S}$ and $H_j|_{\bar S}$ in $\bar S)=($the intersection number of $p^*H_i|_{\bar S}$ and $p^*H_j|_{\bar S}$ in $\tilde S)$. Here the first equality is by (1), the second by generality of $f_i$, so that we can tropicalize the the intersection number, and the last equality is again by generality of $f_i$, so that $H_i$ does not touch $R$. Note that the equalities hold true for $i=j$ as well.

(iii) The Hodge index inequality for the divisors $p^*H_i|_{\bar S},\,i=1,2$, in the smooth compact surface $\tilde S$ translates to the sought Aleksandrov--Fenchel inequality by (2).
\end{proof}

\section{Engineering complete intersections}

\subsection{Sch\"on patchworking} 

Given real open smooth manifolds $U$ and $M\subset U\times\R_{>0}\xrightarrow{t}\R_{>0}$, we can regard the latter as a family $M_{t_0}:=M\cap\{t=t_0\}$ and ask about the topology of a member for $t\to+\infty$ in terms of a compactification of $U\times\R_t$.

We say that $X\supset U\times\R_t$ is a {\it sch\"on compactification}, if

-- the coordinate function $t$ extends to continuous proper $t:X\to(0,+\infty]$;

-- the manifold $U\times\R_t$ is the interior of the smooth manifold $X$ with a normal crossing boundary: each $x\in\partial X$ has a neighborhood $U_x$ with a diffeomorphism $\varphi_x$ to $\R^k_{>0}\times\R^{n-k}$, sending $t$ to a monomial in the standard coordinates of $\R^n$, and $\bar M$ to a coordinate plane.

\begin{rem}
1. In particular, the closure $\bar M\subset X$ is smooth with a normal crossing boundary, and transversal to the strata of $\partial X$.

2. We will call the diffeomorphism
$$\varphi_x:U_x\to\R^k_+\times\R^{n-k}\eqno{(\times)}$$ a {\it monomial chart} near $x$.
\end{rem}

\begin{theor}[The topological patchworking]
If $X$ is a sch\"on compactification of a family $M_t$, there is a homeomorphism of pairs $$(\overline{\{t=t_0\}},\bar M_{t_0})\mbox{ and }(\overline{\{t=+\infty\}},\bar M_{+\infty}),$$
where all closures are taken in $X$, and $t_0$ is large enough.
\end{theor}
\begin{proof}
For a monomial chart $\varphi_x$ of every point $x\in\overline{\{t=+\infty\}}$, reorder the coordinates $(y_1,\ldots,y_n)$ in the target $\R^n$ so that the monomial $t\circ\varphi_x^{-1}$ nontrivially depends on the first $m$ coordinates,  and let $v_x$ be the preimage of the constant vector field $(\underbrace{1,\ldots,1}_m,0,\ldots,0)$.

By construction, the trajectories of this vector field in the neighborhood $U_x$ are transversal to each level $\overline{\{t=t_0\}}$ and tangent to the mainifold $\bar M$ (if the latter is present in $U_x$). 

Since $\overline{\{t=+\infty\}}$ is compact by properness of $t$, we can cover it by finitely many monomial charts $U_{x_i}$, and use a partition of unity to glue the local fields $v_{x_i}$ into a field $v$ in a neighborhood of $\overline{\{t=+\infty\}}$.

By construction, its trajectories are transversal to each $\overline{\{t=t_0\}}$ and tangent to $\bar M$, so the contraction along its trajectories gives the sought homeomorphism.
\end{proof}

A {\it classical setting for patchworking} is when $U=\R^n_{>0}$ is the positive octant, $X$ is the real part of a smooth toric variety $X_P$ associated to a lattice polytope $P\subset\Z^n\times\Z^1$, and the intersections of the manifold $\bar M$ with the strata of $\overline{\{t=+\infty\}}$ (i.e. with repsective orbits of the toric variety) are linear in the following sense.

\begin{defin}
A {\it twisted hyperplane arrangement complement} (THAC) is the preimage of an affine plane in a torus $\CC^m$ under a surjection of tori $\CC^m\to\CC^m$.

A THAC over $\R$ is the intersection of a complex THAC with $\R^m_{>0}$.
\end{defin}

In this classical setting for patchworking, we can construct a piecewise linear {\it patchworking model} for the sought manifold $$\bar M_{t_0}\simeq\bar M_{+\infty}\eqno{(\times\times)}.$$
It will be a piecewsie linear set $Z$  in the {\it upper boundary} $\partial_+ P\subset P$, the topological closed disc defined as the union of all bounded faces of 
$F\in P+\{0\}\times\R_{<0}$.

We shall define it by its intersections with all faces $F\subset\partial_+ P$, arbitrarily ordered by increasing dimension:

-- if the $F$-orbit of the toric variety $X_P$ does not intersect $\bar M$, then we set $Z\cap F:=\varnothing$;

-- if the $F$-orbit of the toric variety $X_P$ intersects $\bar M$, and has smallest possible dimension for this (equal to $\codim M$), then we define $Z\cap F$ as a point in the relative interior of $F$;

-- if the $F$-orbit of the toric variety $X_P$ intersects $\bar M$, and we already know by induction on $\dim F$ the set $Z\cap F'$ for all $F'$ of dimension smaller than $F$, then define $Z\cap F$ as the union of all segments connecting a marked interior point of $F$ with the points of $Z\cap F'$ for $F'\subsetneq F$.

\begin{theor}[Geometric patchworking]\label{geompatch}
In this classical setting for patchworking, the sought manifold $(\times\times)$ is homeomorphic to the patchworking model $Z$ (as defined in the preceding paragraph).
\end{theor}
For the proof, recall that the {\it moment map} of a finite set $A\subset\R^m$, sending $\R^m_{>0}$ to $\conv A$ by
$$y\mapsto \frac{\sum_{a\in A} y^a\cdot a}{\sum_{a\in A} y^a},$$
extends to a homeomorphism of the toric compactification $\bar \R^m_{>0}\subset X_A$ and the convex hull of $A$ (\cite{gkz}).
\begin{proof}
Send $\partial_+P$ to $\overline{\{t=+\infty\}}$ by the moment map $m$ for $A:=($vertices of $P)$, then construct the thought homomorphism $Z\to\bar M_{+\infty}$ on all faces $F$ ordered by increasing dimensions, thanks to the fact that both $Z\cap F$ and $\bar m(M_{+\infty})\cap F$ are discs (the former because it is the intersection of an affine space with an octant, and the sceond by \cite{ardila}). 
\end{proof}

\begin{rem}
1. We do not touch the question of stepping from $\R^n_{>0}$ to $\R^n$, as patchworking the part of the sought variety in each of the $2^n$ octants as above, and then gluing them together, is straightforward.

2. A far more systematic treatment of patchworking for non-complete intersection varieties can be found e.g. in \cite{patchw22} and \cite{patchw23}. We have introduced this subsection just to fix generalities and notation that we will apply to ECIs.
\end{rem}

\subsection{Engineering over $\C$}

We introduce a version of patchworking for ECIs (called {\it engineering}), and use it over both complex numbers (to prove a key lemma for the next section) and real numbers (which allows to construct e.g. real discriminants and reaction networks with prescribed topology).

\begin{defin}
Given a (symbolic) ECI $(\diamondsuit)$, we choose the sign function $s:A\to\{\pm 1\}$, the order function $d:A\to\Q$ in general position, and from these data define the engineering
$$f_i(z,t):=\sum_{a\in A} v_{a,i} \cdot s_a \cdot t^{d_a} \cdot z^a,\eqno{(\diamondsuit\diamondsuit\diamondsuit)}$$
\end{defin}
\begin{utver}
Given a Zariski open $U\subset\C^A$, for a general enough $d$ and large enough $t_0$, the complete intersection $f_1(z,t_0=\cdots=f_k(z,t_0)=0$, regarded as a point in $\C^A$, belongs to $U$.
\end{utver}
\begin{rem}
In particular, for large $t_0$ and general $d$ this complete intersection is sch\"on, because almost all specializations of a given ECI are sch\"on (\cite{e24b}),
\end{rem}
\begin{proof}
This is just a tangled version of the following well known fact: for Zariski open $U\subset\CC^N$, a general 1-parametric subgroup $t^d:\CC^1\to\CC^N$ send all points $t$ except for finitely many to $U$. A sufficient generality for $d$ is to take it outside the tropical fan of the complement of $U$.
\end{proof}
\begin{defin}
For a toric compactification of $\CC^n_z\times\CC^1_t$, its {\it upper orbits} are the orbits, whose (relatively open) cones in the fan of the variety belong to the uper half-space $\R^n\times\R^1_{>0}\subset\R^n\times\R^1$.
\end{defin}
\begin{utver}\label{ECITHAC}
For every engineering $f_1(z,t)=\cdots=f_k(z,t)=0$, there exists a smooth toric compactification $X\supset\CC^n_z\times\CC^1_t$ in which the closure of the set, defined by the engineering equations in $\CC^n_z\times\CC^1_t$, transversally intersects all upper orbits by (reduced, including possibly empty) THACs.
\end{utver}
\begin{proof}
It is enough to prove that, for any linear valuation $l:\Z^n\times\Z\to\Z$ with $l(0,\ldots,0,1)>0$, there exists a cancellation of the system $f_1(z,t)=\cdots=f_k(z,t)=0$ whose leading parts define a reduced THAC. Start constructing such a cancellation $$\tilde f_i(z,t)=f_i(z,t)+g_{i,i-1}f_{i-1}(z,t)+\cdots+g_{i,1}f_1(z,t)$$ from $i=1$ to $k$, choosing constant coefficients $g_{i,j}$ in each step so that the $l$-degree of $\tilde f_i$ is minimal possible. We will notice that, thanks to the general position of the order function $d$, for each $i$ the $v$-leading parts of $\tilde f_1(z,t),\ldots,\tilde f_i(z,t)$ will define THAC, because their support sets will be a sequence of vertex sets of mutually transversal simplices (possibly repeating).
\end{proof}
\begin{rem}\label{remseq}
1. As we see from the proof, every linear valuation $l:\Z^n\times\Z\to\Z$ defines a sequence of transversal simplices that we shall denote by $A^l_i\subset A$.

2. In general, this tuple of simplices depends on the coefficients $v_{a,i}$ of our ECI $(\diamondsuit)$. However, the respective simplices for the ECI $(\tilde f_1,\ldots,\tilde f_k):=(f_1,\ldots,f_k)\cdot($a general complex $k\times k$ matrix$)$ depend only on the respective matroid complete intersection (which is the same for $f$ and $\tilde f$).

Note that at the same time $\{\tilde f=0\}=\{f=0\}$, so we rename $\tilde f$ into $f$ and study it instead of the original ECI throughout this section (c.f. \cite[Section 7]{e24b} for a more direct application of passing from $f$ to $\tilde f$).
\end{rem}

We shall use engineering in the proof of the following equality. Given nondegenerate ECI $f_1=\cdots=f_{k}=0$ in $\CC^n$ with the TCI $(T_i,m_i)$, assume that the closure of $S=\{f_1=\cdots=f_{k-1}=0\}$ in a tropical compactification $X\supset\CC^n$ is smooth (this is the case e.g. if the ECI is sch\"on). Then, for the closure of $C=\{f_1=\cdots=f_{k}=0\}$, its self-intersection number $[\bar C]^{n-k+1}$ in $\bar S$ is defined. 
\begin{sledst}\label{selfint}
In the setting of the paragraph above, the self-intersection number $\bar C^{n-k+1}$ in $\bar S$ equals the tropical $(n-k+1)$-fold self-intersection number of $T_{k}$ in $T_{k-1}$ (defined as $\delta^n(m_1\cdot\ldots\cdot m_{k-1}\cdot m_{k}^{n-k+1})/n!$).
\end{sledst}
\begin{proof}
By proposition \ref{ECITHAC}, an engineering for our ECI is a cohomologically tropical deformation of $(X,S,[\bar C]$ (because every THAC is a twisted HTV). By Theorem \ref{thdef}, this implies the sought equality. 
\end{proof}

\begin{rem}
The general version is: if we extend a nondegenerate ECI $S=\{f_{1,\ldots,k}=0\}$ with the TCI $(T_1,\ldots,T_k, m_1,\ldots,m_k)$ to non-degenerate ECIs $C_j:=\{f_{1,\ldots,k}=f^j_{k+1,\ldots,n_j}=0\}$ with the TCIs $(T_1,\ldots,T_k,T^j_{k+1},\ldots,T^j_{n_j} m_1,\ldots,m_k,m^j_{k+1},\ldots,m^j_{n_j})$, then the intersection number of $\bar C_j$ in $\bar S$ equals $\delta^n(\prod_i m_i\cdot \prod_{i,j} m^j_i)/n!$.
\end{rem}

\subsection{Engineering over $\R$} 

\begin{defin}
1. A {\it real matroid complete intersection} is an oriented matroid $r$ on the ground multiset $A\subset\Z^n$.

2. A real ECI $(\diamondsuit)$ represents the real MCI $(A,r)$, if the vectors $v_a:=(v_{a,i})$ represent $r$.
\end{defin}
\begin{theor}\label{patchreci}
If $f_i(z,t),\,i=1,\ldots,k$ is the engineering of a real ECI $(\diamondsuit)$ defined by the sign function $s$ and the order function $d$, the topological type of $Z_{t_0}:=\{z,|\,f_\bullet(z,t_0)=0\}\subset\R^n_{>0}$ and its closure in any projective sch\"on compactification $X_\Sigma$ is defined by the real MCI of the ECI and the combinatorial data $(s,d,\Sigma)$ for $t_0$ large enough.
\end{theor}
We define the ``zero locus of a real MCI'' and then prove that it is homeomorphic to the sought zero locus of the real ECI.

\vspace{1ex}

Given a real MCI $(A,r),\,A\subset\R^n$ and an order function $d:A\to \Z$ in general position, define its engineering $(\tilde A, r)$ as the matroid $r$ transferred to the ground set $\tilde A:=\{$graph of $d\}\subset\Z^n\times\Z$. Note that $\tilde A$ is a (mono)set, because $d$ is general.

Let $(\tilde T_i,\tilde m_i)$ be the TCI of $(\tilde A, r)$ (as defined in Proposition \ref{proptrop}).

\vspace{1ex}

Notice that, for any linear valuation $l\in T_k$, the projections of $\tilde A\cap\{l=m_i\}\subset\Z^n\times\Z$ to $\Z^n$ form a tuple of transversal simplices $A^l_1,\ldots,A^l_k$. 

\vspace{1ex}

Choose a {\it sch\"on polytope} $P\subset\Z^n\times\Z$ with the following properties:

-- its toric variety $X_P$ is smooth,

-- its dual fan $\tilde\Sigma$ extends  $\Sigma\subset(\Z^n\times\{0\})^*$, and its codimension $k$ skeleton contains $T_k$,

-- for every upper cone $C\in\tilde\Sigma$, all valuations $l$ in the relative interior of $C$ give the same tuple of simplices $A^l_i$, so that. we redenote them by $$A^C_i=\{a^C_{i,0},\ldots,a^C_{i,k_i}\}.$$ 

Note that each simplex $A^C_i$ is normal to $C$ by construction.

We define $\sgn a^C_{i,j}$ as the orientation of the basis $$\bigcup A^C_i\setminus\{a^C_{p,q}\,|\,q=j\mbox{ for }p=i\mbox{ and }0\mbox{ otherwise}\}.$$

Note that, changing $0$ to any other fixed choice in this definition could change all $\sgn a^C_{i,j}$ for given $i$ only simultaneously.

The zero set of our real MCI will live in the {\it upper boundary} $\partial_+ P\subset P$, the topological closed disc defined as the union of all bounded faces of 
$F\in P+\{0\}\times\R_{<0}$.

\vspace{1ex}

Given further the sign function $s:A\to\{\pm 1\}$, 
the {\it zero locus} $Z$ of $(A,r,d,s)$ is defined by its intersections with each face $F$, as follows:

-- if the dual cone of $F$ is not in $\tilde T_k$, the intersection is empty.

-- if the dual cone of $F$ is in $\tilde T_k$, and $F$ has the smallest possible dimension $k$ for this, then we define $Z\cap F$ as an interior point of $F$, whenever, for each $i=1,\ldots,k$, vectors $\sgn(a^C_{i,\bullet})$ and $s(a^C_{i,\bullet})$ equal up to a sign (and set $Z\cap F:=\varnothing$ otherwise).

-- if the dual cone of $F$ is in $\tilde T_k$, and we already know by induction on $\dim F$ the set $Z\cap F'$ for all $F'$ of dimension smaller than $F$, then define $Z\cap F$ as the union of all segments connecting a marked interior point of $F$ with the points of $Z\cap F'$ for $F'\subsetneq F$.

We can now prove a constructive version of Theorem \ref{patchreci}.
\begin{utver}
In the setting of Theorem \ref{patchreci}, the real engineering $\bar Z_{t_0}\subset X_\Sigma$ is homeomorphic to the zero locus $Z$ of the combinatorial data $(A,r,d,s)$ (as defined in the preceding paragraph).
\end{utver}
\begin{proof}
This is the special case of Theorem \ref{geompatch}: one just needs to compare the definition of the matroid zero locus $Z$ above and the patchworking model $Z$ defined in  Theorem \ref{geompatch}. The only non-tautological (but straightforward) point in this comparison is that the real THAC in the $F$-orbit is non-empty, iff the face $F$ satisfies $\sgn(a^C_{i,\bullet})=\pm s(a^C_{i,\bullet}).$
\end{proof}

\begin{rem}
The zero locus $Z$ of a real MCI is a topological manifold with a boundary: if the MCI can be represented by a real ECI, this follows from the proposition; but even if we start with a non-representable oriented matroid $r$, this follows by \cite{ardila}.
\end{rem}

\section{Combinatorics of reducibility}

\subsection{Towards a criterion} We give an inductive (by the number of equations) way to check whether an ECI is irreducible.
\begin{theor}\label{thirr}
A complex sch\"on ECI $f_1=\cdots=f_k=0$ with the TCI $(T_i,m_i)$ is irreducible, if $f_1=\cdots=f_{k-1}=0$ is irreducible, and the self-intersection of $T_k$ in $T_{k-1}$ (defined as $\delta^2(m_k^2\cdot T_{k-1})$) is positive in one of the fallowing increasingly mild but involved senses:

i) $m_k\geqslant 0$ and not identical 0 on the fan $T_{k-1}$;

ii) all weights of the fan $\delta^2(m_k^2\cdot T_{k-1})$ are non-negative, and it is not 0;

iii) this fan has a positive intersection number with the tropical fan of some other ECI $G$ (of complimentary dimension).
\end{theor}
There is a weaker criterion in \cite{e24a} for an arbitrary sch\"on SCI, and an extension \cite{zhizhin} of the irreducibility from the classical Newton polytope theory \cite{kh16}, given in terms of the MCI. Our theorem does not cover it, but is more realistic to check in terms of the TCI.

For the proof, we need the following observation.
\begin{lemma}
If a sch\"on ECI curve $f_1=\cdots=f_{n-1}=0$ in $\CC^n$ has several components $C_i$, then 

1. The fan $\Trop C_i$ does not depend on $i$;

2. The self-intersection number $\bar C_i\circ \bar C_i$ in the sch\"on compactification $S$ of the surface $f_1=\cdots=f_{n-2}=0$ does not depend on $i$;

3. $\bar C_i\circ \bar C_j=0$ for $i\ne j$.
\end{lemma}
\begin{proof}
3. By sch\"onness, $\bar C_i$ and $\bar C_j$ are smooth non-intersecting curves in the compact surface $S$.

1 and 2. In the space $\C^A$: each point $c\in\C^A$ can be identified with a specialization of the given symbolic ECI
$$f_i(z):=\sum_{a\in A} v_{a,i} c_a z^a, i=1,\ldots,k-1.$$ Consider the tautological space $V:=\{(f_\bullet,z)\,|\,f_\bullet(z)=0\}\subset\C^A\times\CC^n$. It is a vector bundle over $\CC^n$, thus irreducible.

Pick a Zariski open set of sch\"on nondegenerate realizations $U\subset\C^A$, denote its preimage in $V$ by $W$, and notice that:

i. the set $W$ is connected, because $V$ is irreducible;

ii. the projection $W\to U$ is a locally trivial fibration by \cite{e24a}[Proposition 2.4].

We wish to prove that the discrete invariants of the components $C_1$ and $C_2$ of a given sch\"on nondegenerate specialization $f_{\bullet,0}=0$ coincide. Pick a point $z_i\in C_i$, connect the respective points $(f_{\bullet,0},z_i)$ with a path in $W$ (using (i)), and denote its projection in $U$ (which is a loop from $f_0$) by $\gamma$. Along this loop, the discrete topological invariants of components of the ECI does not change (by (ii)), but $C_1\ni z_1$ travels to the component containing $z_2$, which is $C_2$.
\end{proof}
{\it Proof of the theorem.} Intersecting with $G$ from (iii), we reduce the general case to $k=n-1$.

In this case, adopting the notation from the lemma, and assuming towards the contradiction that we have more than one component $C_1,\ldots,C_k$, we notice $0<\bar C\circ\bar C=\bar \sum C_i\circ\bar C_i=k\, \bar C_1\circ\bar C_1$. This allows to use the Hodge index inequality $0=(\bar C_1\circ\bar C_2)^2\geqslant (\bar C_1\circ\bar C_1)\cdot (\bar C_2\circ\bar C_2)=(\bar C_1\circ\bar C_1)^2$, leading to a contradiction. $\hfill\square$

\subsection{Towards a classification}

We are interested to classify TCIs that can be a tropicalization of a reducible ECI, similarly to the classification in the classical Newton polytope theory \cite{kh16}.
While the complete answer seems out of reach, we will try to have an idea of what can be done in this direction.
Specifically, we will try to classify  pairs of tropical fans $(T,F)$ of the following type, with a view towards having $T=T_{n-2}$ and $F=\delta^2(m_{n-1}^2)$ in the setting of Theorem \ref{thirr}. 
\begin{defin}
A pair of tropical fans (with non-negative weights) $(T,F)$ in $\R^n$ is said to be {\it interesting}, if $T\circ F=0$, and $T\setminus\{0\}$ connected.
\end{defin}

To classify interesting pairs $(T,F)$ in $\R^n$, let $[T]$ be the union of planes generated by the maximal dimensional cones of $T$, and the same for $[F]$. Then the pair of sets $({\mathbb P}[F],{\mathbb P}[T])$ in $\RP^{n-1}$ is interesting in the following sense.

\begin{defin}
Finite sets of p-dimensional planes $L_i$ and 1-dimensional planes $R_j$ in $\RP^{p+q+1}$ form an {\it interesting pair}, if each $L_i$ intersects each $R_j$.
\end{defin} 

\begin{observ}
1. The classification of interesting pairs of fans reduces to the classification of the interesting pairs  $(L_\bullet,R_\bullet)$ in $\RP^{p+q+1}$ with $\bigcup L$ connected.

2. The latter looks achievable.
\end{observ}
We do it in the simplest case $p=q=1$.

\begin{utver}
Every interesting pair $(L_\bullet,R_\bullet)$ in $\RP^3$ with $\bigcup L$ connected falls into one of the following types (up to interchanging $L$ and $R$).

1. All lines $L_i$ and $R_j$ pass through the same point $a$.

2. All lines $L_i$ and $R_j$ belong to the same plane $P$.

3. All lines $L_i$ pass through the same point $a$ and belong to the same plane $a\in P$, each line $R_j$ either belongs to $P$ or contains $a$.

4. For planes $P_1$ and $P_2$ and points $a_1$ and $a_2$ in their intersection, 

-- each line $L_i$ belongs to $P_1$ and contains $a_1$, or belongs to $P_2$ and contains $a_2$;

-- each line $R_i$ belongs to $P_1$ and contains $a_2$, or belongs to $P_2$ and contains $a_1$.

5. The family $L_\bullet$ consists of one line.

6. The family $L_\bullet$ consists of two skew lines $L-1$ and $L-2$, and each line $R_j$ intersects both.
\end{utver}
\begin{proof}
We easily exhaust the cases when all lines $R_j$ contain the same point or all belong to the same plane (arriving at the cases (1-4)). So we exclude them from consideration (and the same for $L_i$).

The only remaining possibility is that there are two skew lines among them, say $R_1$ and $R_2$, intersecting each $L_i$. 
Since $\bigcup L_i$ is connected and not in a plane, it contains lines $AB$, $BC$ and $CD$, with $A,C\in R_1$ and $B,D\in R_2$.

On the other hand, there is one more line $R_3$ (otherwise we arrive at case (6)), and since it intersects the lines $AB$, $BC$ and $CD$, it is contained in the plane $ABC$ or $BCD$. This falls into case (4) with $a_1$ and $a_2$ being $B$ and $C$, and $P_1$ and $P_2$ being $ABC$ and $BCD$.
\end{proof}
\begin{sledst}
Every interesting pair of fans $(F,T)$ in $\R^4$ falls into one of the following types (up to permuting $F$ and $T$):

1. Both fans belong to a 3-plane.

2. Both fans pull back from a projection $\R^4\to\R^3$.

3. The fan $T$ is a pull back of a fan in a 2-plane $P\subset\R^3$ under a projection $p:\R^4\to\R^3$. The fan $F$ is the sum of a fan in $p^{-1}(P)$ and the pull back of a fan in $\R^3$.

4. There are 3-planes $P_1$ and $P_2$ intersecting by a span of lines $L_1$ and $L_2$ such that 

-- the fan $F$ is the sum of pull backs of a fan in $P_1/L_1$ and a fan in $P_2/L_2$;

-- the fan $T$ is the sum of pull backs of a fan in $P_1/L_2$ and a fan in $P_2/L_1$.

5. The fan $F$ is a 2-plane $P$ (with some weight), and each cone of $T$ has at least one boundary ray in $P$.

6. The fan $F$ consists of two transversal 2-planes $P_1$ and $P_2$, and each cone of $T$ has one boundary ray in $P_1$ and the other one in $P_2$.
\end{sledst}

\vspace{1cm}
\noindent
London Institute for Mathematical Sciences, UK \\
\textit{Email}: aes@lims.ac.uk

\end{document}